\documentclass[12pt]{article}
\usepackage[utf8]{inputenc}
\usepackage{amsmath}
\usepackage{amsthm}
\usepackage{amssymb}
\usepackage[affil-it]{authblk}
\newcommand{\Z}{\mathbb Z}
\newcommand{\N}{\mathbb N}

\newcommand{\ZM}{\mathbb Z[\sqrt m]}
\newcommand{\NM}{\mathbb N[\sqrt m]}
\newtheoremstyle{theoremstyle}{}{}{\itshape}{}{\scshape}{.}{ }{\textbf{#1\ #2}}  
\theoremstyle{theoremstyle}
\newtheorem{theorem}{Theorem}
\newtheorem{lem}{Lemma}

\newtheorem{cor}{Corollary}

\title{2-Variable Frobenius Problem in \(\Z[\sqrt M]\)}
\author{\LARGE{Doyon Kim}\footnote{AMS(MOS) Subject Classification: 11D07, 11R11, 16B99}\footnote{This work was supported by NSF grant no. 1262930, and was completed during and after the 2015 summer Research Experience for Undergraduates in Algebra and Discrete Mathematics at Auburn University.}
\\ \large{\textnormal{Auburn University}}
\\ \large\textnormal{dzk0028@auburn.edu}}
\date{July 2015}
\begin{document}

\maketitle
\begin{abstract}
Suppose that \(m\) is a positive integer, not a perfect square. We present a formula solution to the 2-variable Frobenius problem in \(\ZM\) of the "first kind" ([3]).
\end{abstract}

\section{Introduction}
\ \ \ \ Let \(\Z\) denote the set of integers, and \(\N\) the set of non-negative integers; \(z_1, \dots, z_n\in\Z\) are relatively prime, or coprime, if they have no common divisor in \(\Z\) other than \(\pm 1\). Given \(a_1, \dots, a_n\in\N\setminus\{0\}\), let \(SG(a_1,\dots,a_n)=\{\lambda_1a_1+\cdots+\lambda_na_n|\lambda_1,\dots,\lambda_n\in\N\}\); \(SG(a_1,\dots,a_n)\) is the \textit{semigroup} generated by \(a_1, \dots, a_n\). \par
If \(a_1, \dots , a_n\in\N\setminus\{0\}\) are coprime, then \(SG(a_1,\dots,a_n)\) contains a tail of \(\N\), \(\{g,g+1,\dots\}=g+\N\). A Frobenius problem in \(\N\) is the following: \textit{For coprime positive integers \(a_1,\dots,a_n\), find the smallest \(g=g(a_1,\dots,a_n)\) such that \(g+\N\subseteq SG(a_1,\dots,a_n)\).} These are called Frobenius problems because Frobenius provided a beautiful proof of a formula solution in the case \(n=2\): \(g(a_1,a_2)=(a_1-1)(a_2-1)\). No such formula has been found for the cases \(n\geq 3\). \par
In [1], [2] and [3], Frobenius problems in different rings are proposed. In the original Frobenius problem in \(\Z\), \(a_1, \dots, a_n\) are chosen from \(\N\), and the necessary and sufficient condition of \(a_1, \dots , a_n\in\N\setminus\{0\}\) to have a solution of Frobenius problem, \(g=g(a_1,\dots,a_n)\) such that \(g+\N\subseteq SG(a_1,\dots,a_n)\), is that they are coprime. Finding such \(g\) is equivalent to describing the set \(Frob(a_1,\dots,a_n)=\{w\in\Z\mid w+\N\in SG(a_1,\dots,a_n)\}\), because once we find \(g\) we automatically know that \(Frob(a_1,\dots,a_n)=g+\N\). Solving Frobenius problems in a different ring is similar to solving these problems in \(\Z\), but there are few more things to check to ensure that the Frobenius problem has a solution. 
\par
\section{Frobenius problems in \(\ZM\)}
\ \ \ \ Suppose \(R\) is a commutative ring with multiplicative identity \(1\). We say that a sequence \((\alpha_1,\dots,\alpha_n)\in R^n\) spans \(1\) in \(R\) if and only if for some \(\lambda_1,\dots,\lambda_n\in R\), \(\lambda_1\alpha_1+\cdots+\lambda_n\alpha_n=1\). (In other words, the ideal generated by \(\alpha_1,\dots,\alpha_n\) in \(R\) is \(R\).) \par
Let \(m\) be a positive integer which is not a perfect square, and \(\ZM=\{a+b\sqrt m\mid a,b\in\Z\}\). Because \(\sqrt m\) is irrational, for each \(\alpha\in\ZM\) there exist unique \(a,b\in\Z\) such that \(\alpha=a+b\sqrt m\). Let \(\ZM^{+}=\ZM\cap [0,\infty)\) and \(\NM=\{a+b\sqrt m\mid a,b\in\N\}\). For \(\alpha_1,\dots,\alpha_n\in\ZM\), let \(SG(\alpha_1,\dots,\alpha_n)=\{\lambda_1\alpha_1+\cdots+\lambda_n\alpha_n\mid \lambda_1,\dots,\lambda_n\in\NM\}\) and let \(SG'(\alpha_1,\dots,\alpha_n)=\{\lambda_1\alpha_1+\cdots+\lambda_n\alpha_n\mid \lambda_1,\dots,\lambda_n\in\ZM^{+}\}\). Also, let \(Frob(\alpha_1,\dots,\alpha_n)=\{w\in\ZM\mid w+\NM\subseteq SG(\alpha_1,\dots,\alpha_n)\}\) and let \(Frob'(\alpha_1,\dots,\alpha_n)=\{w\in\ZM\mid w+\ZM^{+}\subseteq SG'(\alpha_1,\dots,\alpha_n)\}\). There are two kinds of Frobenius problem in \(\ZM\). One kind is, for \(\alpha_1,\dots,\alpha_n\in\NM\setminus\{0\}\), to describe the set \(Frob(\alpha_1,\dots,\alpha_n)\). The other kind is, for \(\beta_1,\dots,\beta_n\in\ZM^{+}\), to describe the set \(Frob'(\beta_1,\dots,\beta_n)\). \par
In [3] it is shown, and is easy to see directly, that in each kind of Frobenius problem in \(\ZM\), for \(Frob(\alpha_1,\dots,\alpha_n)\neq\emptyset\) (or \(Frob'(\alpha_1,\dots,\alpha_n)\neq\emptyset\)) it is necessary that the sequence \((\alpha_1,\dots,\alpha_n)\) span \(1\) in \(\ZM\). In [2] it is proved that given \(\beta_1\), \dots, \(\beta_n\in\ZM^{+}\), \(Frob'(\beta_1,\dots,\beta_n)\) is not empty if \((\beta_1,\dots,\beta_n)\) spans \(1\) in \(\ZM\), and in [3] it was shown that if \(\beta_1\), \dots, \(\beta_n\in\ZM^{+}\) spans \(1\) then \(Frob'(\beta_1,\dots,\beta_n)=\ZM^{+}\). So the second kind of the problem, which is on \(\ZM^{+}\), is pretty much solved. \par
We will focus on the first kind of Frobenius problem, with \(n=2\). Let \(\alpha_1,\alpha_2\in\NM\setminus\{0\}\), and \(SG(\alpha_1,\alpha_2)=\{\lambda_1\alpha_1+\lambda_2\alpha_2\mid \lambda_1,\lambda_2\in\NM\}\). Our goal is to describe the set \(Frob(\alpha_1,\alpha_2)=\{w\in\ZM\mid w+\NM\subseteq SG(\alpha_1,\alpha_2)\}\). \par 
In [2] it is proved that if \(\alpha_1,\dots,\alpha_n\in\NM\setminus\{0\}\) span 1 in \(\ZM\), then \(Frob(\alpha_1,\dots,\alpha_n)\neq\emptyset\) if and only if some \(\alpha_i\) has either rational or irrational part \(0\), and [3] presented the solutions for the cases where every \(\alpha_1,\dots,\alpha_n\) has either rational or irrational part \(0\). Thanks to [3], we already know that for \(m\in \N\setminus\{n^2|n\in \N\}\) and \(a,b\in \N\setminus\{0\}\), if \(a\) and \(b\) are coprime then \(Frob(a,b)=(a-1)(b-1)(1+\sqrt m)+\NM,\) and if \(a\) and \(bm\) are coprime then \(Frob(a,b\sqrt m)=(a-1)(b\sqrt m-1)(1+\sqrt m)+\NM\); \(Frob(a\sqrt m,b\sqrt m)\) is always empty because \(a\sqrt m,b\sqrt m\) cannot span \(1\). \par 
So the two remaining cases are \((\alpha_1,\alpha_2)=(a, b+c\sqrt m)\) and \((\alpha_1,\alpha_2)=(a\sqrt m, b+c\sqrt m)\), where \(a,b\) and \(c\) are positive integers. For each case, we will first find the conditions for such \((\alpha_1,\alpha_2)\) to span \(1\) and will find \(Frob(\alpha_1,\alpha_2)\). This will complete the solutions of the first Frobenius problem in \(\ZM\) in the case \(n=2\). The answer, by the way, is that if \(\alpha,\beta\in \NM\setminus\{0\}\) span \(1\), then \(Frob(\alpha,\beta)=(\alpha-1)(\beta-1)(\sqrt m+1)+\NM\), which agrees with the results in [3] in the special cases solved there.
\section{\((a,b+c\sqrt m)\)}

Let \(a,b,c\in \N\setminus\{0\}\) and let \(m\in \N\setminus\{n^2|n\in \N\}\).
\begin{lem} %lemma 1
\((a, b+c\sqrt m)\) spans 1 in \(\Z[\sqrt m]\) if and only if \((a,b^2-c^2m)\) spans 1 in \(\Z\).
\begin{proof}
Suppose \((a,b^2-c^2m)\) spans \(1\) in \(\Z\). Then \[a\lambda+(b^2-c^2m)\mu=1\] for some \(\lambda, \mu\in \Z\). Let \[x=\lambda+0\cdot\sqrt m=\lambda\in \Z[m],\quad y=(b-c\sqrt m)\mu\in \Z[m].\] Then \[ax+(b+c\sqrt m)y=a\lambda+(b+c\sqrt m)(b-c\sqrt m)\mu=a\lambda+(b^2-c^2m)\mu=1.\]
Now, suppose \((a, b+c\sqrt m)\) spans \(1\) in \(\Z[\sqrt m]\). Then \[a(x+y\sqrt m)+(b+c\sqrt m)(z+w\sqrt m)=1\] for some \(x+y\sqrt m,z+w\sqrt m\in \Z[\sqrt m]\). Such \(x,y,z,w\in \Z\) satisfies 
\[ax+bz+cmw=1,\]
\[ay+cz+bw=0.\]
Let \(\gcd(a,b^2-c^2m)=g\). We have
\[
c(ax+bz+cmw)-b(ay+cz+bw)=a(cx-by)-(b^2-c^2m)w=c
\]
so \(g\mid c\), and, 
\[
b(ax+bz+cmw)-cm(ay+cz+bw)=a(bx-cmy)+(b^2-c^2m)z=b
\]
so \(g\mid b\). Since \(g\mid a\), \(g\mid b\) and \(g\mid c\), \(g\mid ax+bz+cwm=1\). So \(g=1\), and therefore \((a,b^2-c^2m)\) spans 1 in \(\Z\).
\end{proof}
\end{lem}

\begin{lem} %lemma 2
For \(a,b,c \in \N\setminus\{0\}\), \(m\in \N\setminus\{n^2|n\in \Z\}\) and \(\gcd(a,b^2-c^2m)=1\), \[ax+(b+c\sqrt m)y=A+B\sqrt m\] has a solution in \(\ZM\) for every \(A+B\sqrt m\in \Z[\sqrt m]\).
\begin{proof}
Let \(A+B\sqrt m\in \Z[\sqrt m]\). By Lemma 1, \[ax+(b+c\sqrt m)y=1\] for some \(x, y\in \ZM\). Then \(x(A+B\sqrt m)\in\ZM\), \(y(A+B\sqrt m)\in \ZM\) and they satisfy \[ax(A+B\sqrt m)+(b+c\sqrt m)y(A+B\sqrt m)=A+B\sqrt m.\]
\end{proof}
\end{lem}

From now on, we assume that \(a\), \(b\), \(c\) are positive integers, \(m\) is a positive integer that is not a perfect square, \(\gcd(a,b^2-c^2m)=1\), and \(A,B\in \Z\).
\begin{lem}
If \((x_0,y_0, z_0, w_0)\in \Z^4\) is a solution of \[a(x+y\sqrt m)+(b+c\sqrt m)(z+w\sqrt m)=A+B\sqrt m,\] then every other solution \((x',y',z',w')\in \Z^4\) satisfies \[a\mid z_0-z',\quad a\mid w_0-w'.\]
\begin{proof}
Let \((x_0,y_0, z_0, w_0)\in \Z^4\) be a solution of \[a(x+y\sqrt m)+(b+c\sqrt m)(z+w\sqrt m)=A+B\sqrt m.\] Then \(x_0\), \(y_0\), \(z_0\), \(w_0\in\Z\) satisfies \[ax_0+bz_0+cmw_0=A,\] \[ay_0+cz_0+bw_0=B.\]
Let \((x',y', z', w')\in \Z^4\) be another solution. \(x'\), \(y'\), \(z'\), \(w'\in\Z\) also satisfies \[ax'+bz'+cmw'=A,\] \[ay'+cz'+bw'=B.\]
So we have \[a(x_0-x')+b(z_0-z')+cm(w_0-w')=0,\] \[a(y_0-y')+c(z_0-z')+b(w_0-w')=0.\]
Since
\[ab(x_0-x')=-b^2(z_0-z')-bcm(w_0-w') \]
and
\[acm(y_0-y')=-c^2m(z_0-z')-bcm(w_0-w'), \]
we have
\[a(b(x_0-x')-cm(y_0-y'))=-(b^2-c^2m)(z_0-z').\]
Since \(gcd(a,b^2-c^2m)=1\), \(a\mid z_0-z'\). Also, we have
\[ac(x_0-x')=-bc(z_0-z')-c^2m(w_0-w') \]
and
\[ab(y_0-y')=-bc(z_0-z')-b^2(w_0-w'), \]
so \[a(c(x_0-x')-b(y_0-y'))=(b^2-c^2m)(w_0-w').\]
Since \(gcd(a,b^2-c^2m)=1\), \(a\mid w_0-w'\).
\end{proof}
\end{lem}
\begin{lem}
If \((x_0,y_0, z_0, w_0)\in \Z^4\) is a solution of \[a(x+y\sqrt m)+(b+c\sqrt m)(z+w\sqrt m)=A+B\sqrt m,\] then for every \(k\in \Z\), \[x'=x_0-(b-cm)k,\quad y'=y_0+(b-c)k,\quad z'=z_0+ak,\quad w'=w_0-ak\] is also a solution.
\begin{proof}
The proof is straightforward.
\end{proof}
\end{lem}
\begin{lem}
If \((x_0,y_0, z_0, w_0)\in \Z^4\) is a solution of \[a(x+y\sqrt m)+(b+c\sqrt m)(z+w\sqrt m)=A+B\sqrt m,\] then for every \(l\in \Z\), \[x'=x_0+bl,\quad y'=y_0+cl,\quad z'=z_0-al,\quad w'=w_0\] is also a solution.
\begin{proof}
The proof is straightforward.
\end{proof}
\end{lem}
\begin{cor}
There is a unique solution \((\bar x,\bar y,\bar z, \bar w)\) of \[a(x+y\sqrt m)+(b+c\sqrt m)(z+w\sqrt m)=A+B\sqrt m\] with \(0\leq \bar z, \bar w<a\).
\begin{proof}
Let  \((x_0,y_0, z_0, w_0)\) be a solution. Let \(k=\lfloor \frac {w_0}{a}\rfloor\). Then \(k\in \Z\), and  \[x'=x_0-(b-cm)k,\quad y'=y_0+(b-c)k,\quad z'=z_0+ak,\quad w'=w_0-ak\] is a solution. By the choice of \(k\), \(0\leq w'<a\) and by Lemma 3, such \(w'\) is unique. \\
Now, let \(l=\lfloor \frac {z'}{a}\rfloor\). Then \(l\in \Z\), and \[\bar x=x'+bl,\quad \bar y=y'+cl,\quad \bar z=z'-al,\quad \bar w=w'\] is also a solution, and by the choice of \(l\), \(0\leq \bar z<a\). By Lemma 3, such \(\bar z\) is unique. So \((\bar x,\bar y,\bar z, \bar w)\) is a solution with \(0\leq \bar z, \bar w<a\). This is a unique such solution, because if \(\bar z\) and \(\bar w\) are fixed, so are \(\bar x\) and \(\bar y\).
\end{proof}
\end{cor}
\begin{lem}
\(a(x+y\sqrt m)+(b+c\sqrt m)(z+w\sqrt m)=A+B\sqrt m\) has a solution in \(\NM\) if and only if the unique solution \((\bar x,\bar y,\bar z, \bar w)\) with \(0\leq \bar z, \bar w<a\) satisfies \(0\leq \bar x\), \(0\leq \bar y\).
\begin{proof}
Suppose \(\bar x, \bar y\geq 0\). Then \(\bar x+\bar y\sqrt m\), \(\bar z+\bar w\sqrt m\) is a solution in \(\NM\). Now, suppose either \(\bar x<0\) or \(\bar y<0\).  Let \((x', y', z', w')\) be another solution. If \(z'<\bar z\) then \(z'<0\), and if \(w'<\bar w\) then \(w'<0\), by Lemma 3. If \(z'\geq \bar z\) and \(w'\geq \bar w\), then \[x'=\frac{1}{a}(A-bz'-cmw')\leq \frac{1}{a}(A-b\bar z-cm\bar w)=\bar x,\quad \textrm{and}\]
\[y'=\frac{1}{a}(B-cz'-bw')\leq \frac{1}{a}(B-c\bar z-b\bar w)=\bar y\]
so either \(x'<0\) or \(y'<0\).
\end{proof}
\end{lem}
\begin{theorem}
If \(A\geq (a-1)(b-1+cm)\) and \(B\geq (a-1)(b-1+c)\), then \(a(x+y\sqrt m)+(b+c\sqrt m)(z+w\sqrt m)=A+B\sqrt m\) has a solution in \(\NM\).
\begin{proof}
Let \(A\geq (a-1)(b-1+cm)\) and \(B\geq (a-1)(b-1+c)\). Consider the equation \(a(x+y\sqrt m)+(b+c\sqrt m)(z+w\sqrt m)=A+B\sqrt m\). There is a unique solution \((\bar x,\bar y,\bar z, \bar w)\) with \(0\leq \bar z, \bar w<a\); \((\bar x,\bar y,\bar z, \bar w)\) satisfies 
\[a\bar x+b\bar z+cm\bar w=A,\]
\[a\bar y+c\bar z+b\bar w=B.\] 
Since \(0\leq \bar z, \bar w\leq a-1\), we have 
\[a\bar x=A-b\bar z-cm\bar w\geq (a-1)(b-1+cm)-b(a-1)-cm(a-1)=1-a\] so
\[\bar x\geq -1+\frac{1}{a}>-1 \iff \bar x\geq 0.\] Likewise, 
\[a\bar y=B-c\bar z-b\bar w\geq (a-1)(b-1+c)-c(a-1)-b(a-1)=1-a\] so
\[\bar y\geq -1+\frac{1}{a}>-1 \iff \bar y\geq 0.\] Therefore \(\bar x+\bar y\sqrt m\in \NM\), \(\bar z+\bar w\sqrt m\in \NM\) and \(\bar x+\bar y\sqrt m\), \(\bar z+\bar w\sqrt m\) is a solution of \(a(x+y\sqrt m)+(b+c\sqrt m)(z+w\sqrt m)=A+B\sqrt m\) in \(\NM\).
\end{proof}
\end{theorem}
\begin{theorem}
\(Frob(a,b+c\sqrt m)=(a-1)(b+c\sqrt m-1)(1+\sqrt m)+\NM\).
\begin{proof}
Note that \[(a-1)(b+c\sqrt m-1)(\sqrt m+1)=(a-1)(b-1+cm)+(a-1)(b-1+c)\sqrt m.\]
In view of Theorem 1, to show that \[Frob(a,b+c\sqrt m)=(a-1)(b+c\sqrt m-1)(\sqrt m+1)+\NM\] it suffices to show that there is no \(\beta\in \Z\) such that \[(a-1)(b-1+cm)-1+\beta\sqrt m\in Frob(a,b+c\sqrt m)\] and that there is no \(\alpha\in \Z\) such that \[\alpha+((a-1)(b-1+c)-1)\sqrt m\in Frob(a,b+c\sqrt m).\]
Let \[A=(a-1)(b-1+cm)-1=ab-a-b+acm-cm.\] 
Take arbitrary \(k\in \N\), and let \[B=(a-1)(b+c)+ak.\] Then \((-1,k,a-1,a-1)\) is the unique solution \((\bar x, \bar y, \bar z, \bar w)\) of \[a(x+y\sqrt m)+(b+c\sqrt m)(z+w\sqrt m)=A+B\sqrt m\] with \(0\leq \bar z, \bar w<a\). By Lemma 6, the equation does not have a solution in \(\NM\). Since the choice of \(k\) was arbitrary, \(B\) can be arbitrarily large. Therefore there is no \(\beta\in \Z\) such that \[(a-1)(b-1+cm)-1+\beta\sqrt m\in Frob(a,b+c\sqrt m).\]
Now, let \[B=(a-1)(b-1+c)-1.\]
Take arbitrary \(l\in \N\), and let \[A=(a-1)(b+cm)+al.\]
Then \((l,-1,a-1,a-1)\) is the unique solution \((\bar x, \bar y, \bar z, \bar w)\) of \[a(x+y\sqrt m)+(b+c\sqrt m)(z+w\sqrt m)=A+B\sqrt m\] with \(0\leq \bar z, \bar w<a\). By Lemma 6, the equation does not have a solution in \(\NM\). Since the choice of \(l\) was arbitrary, \(A\) can be arbitrarily large. Therefore there is no \(\alpha\in \Z\) such that \[\alpha+((a-1)(b-1+c)-1)\sqrt m\in Frob(a,b+c\sqrt m).\]
\end{proof}
\end{theorem}
\section{\((a\sqrt m,b+c\sqrt m)\)}
Let \(a,b,c\in \N\setminus\{0\}\) and let \(m\in \N\setminus\{n^2|n\in \Z\}\).
\begin{lem} %lemma 7
\((a\sqrt m, b+c\sqrt m)\) spans 1 in \(\Z[\sqrt m]\) if and only if \((am,b^2-c^2m)\) spans 1 in \(\Z\).
\begin{proof}
Suppose \((am,b^2-c^2m)\) spans \(1\) in \(\Z\). Then \[am\lambda+(b^2-c^2m)\mu=1\] for some \(\lambda, \mu\in \Z\). Let \[x=\sqrt m\lambda\in \ZM,\quad y=(b-c\sqrt m)\mu\in \ZM.\] Then \[a\sqrt mx+(b+c\sqrt m)y=1.\]
Now, suppose \((a\sqrt m, b+c\sqrt m)\) spans \(1\) in \(\Z[\sqrt m]\). Then \[a\sqrt m(x+y\sqrt m)+(b+c\sqrt m)(z+w\sqrt m)=1\] for some \(x+y\sqrt m,z+w\sqrt m\in \Z[\sqrt m]\); then \(x,y,z,w\in \Z\) satisfy
\[amy+bz+cmw=1,\]
\[ax+cz+bw=0.\]
Then \((ay+cw)m+zb=1\) so \(\gcd(m,b^2-c^2m)=\gcd(m,b)=1\).
Let \(\gcd(a,b^2-c^2m)=g\). We have
\[
c(amy+bz+cmw)-b(ax+cz+bw)=a(cmy-bx)-(b^2-c^2m)w=c
\]
so \(g\mid c\), and, 
\[
b(amy+bz+cmw)-cm(ax+cz+bw)=a(bmy-cmx)+(b^2-c^2m)z=b
\]
so \(g\mid b\). Since \(g\mid a\), \(g\mid b\) and \(g\mid c\), \(g\mid ax+bz+cwm=1\). Therefore \(g=1\). Since \(gcd(a,b^2-c^2m)=1\) and \(\gcd(m,b^2-c^2m)=1\), we conclude that \(\gcd(am,b^2-c^2m)=1\). Therefore \((am,b^2-c^2m)\) spans 1.
\end{proof}
\end{lem}
\begin{lem} %lemma 8
For \(a,b,c \in \N/\{0\}\), \(m\in \N/\{n^2|n\in \Z\}\) and \(\gcd(am,b^2-c^2m)=1\), \[a\sqrt mx+(b+c\sqrt m)y=A+B\sqrt m\] has a solution in \(\ZM\) for every \(A+B\sqrt m\in \Z[\sqrt m]\).
\begin{proof}
This follows directly from Lemma 7. We skip the proof, as it is almost exactly same as the proof of Lemma 2.
\end{proof}
\end{lem} 
From now on, we assume that \(a\), \(b\), \(c\) are positive integers, \(m\) is a positive integer that is not a perfect square, \(\gcd(am,b^2-c^2m)=1\), and \(A,B\in \Z\).

\begin{lem} %lemma 9
If \((x_0,y_0, z_0, w_0)\in \Z^4\) is a solution of \[a\sqrt m(x+y\sqrt m)+(b+c\sqrt m)(z+w\sqrt m)=A+B\sqrt m,\] then every other solution \((x',y',z',w')\in \Z^4\) satisfies \[am\mid z_0-z',\quad a\mid w_0-w'.\]
\begin{proof}
Let \((x_0,y_0, z_0, w_0)\in \Z^4\) be a solution of \[a\sqrt m(x+y\sqrt m)+(b+c\sqrt m)(z+w\sqrt m)=A+B\sqrt m.\] Then \(x_0\), \(y_0\), \(z_0\), \(w_0\) satisfies \[amy_0+bz_0+cmw_0=A,\] \[ax_0+cz_0+bw_0=B.\]
Let \((x',y', z', w')\in \Z^4\) be another solution. \(x'\), \(y'\), \(z'\), \(w'\) satisfies \[amy'+bz'+cmw'=A,\] \[ax'+cz'+bw'=B.\]
So we have \[am(y_0-y')+b(z_0-z')+cm(w_0-w')=0,\] \[a(x_0-x')+c(z_0-z')+b(w_0-w')=0.\]
Since
\[abm(y_0-y')=-b^2(z_0-z')-bcm(w_0-w') \]
and
\[acm(x_0-x')=-c^2m(z_0-z')-bcm(w_0-w'), \]
we get
\[am(b(y_0-y')-c(x_0-x'))=-(b^2-c^2m)(z_0-z').\]
Since \(gcd(am,b^2-c^2m)=1\), \(am\mid z_0-z'\). Also, we have
\[acm(y_0-y')=-bc(z_0-z')-c^2m(w_0-w') \]
and
\[ab(x_0-x')=-bc(z_0-z')-b^2(w_0-w'), \]
so \[a(cm(y_0-y')-b(x_0-x'))=(b^2-c^2m)(w_0-w').\]
Since \(gcd(a,b^2-c^2m)=1\), \(a\mid w_0-w'\).
\end{proof}
\end{lem}
\begin{lem} %lemma 10
If \((x_0,y_0, z_0, w_0)\in \Z^4\) is a solution of \[a\sqrt m(x+y\sqrt m)+(b+c\sqrt m)(z+w\sqrt m)=A+B\sqrt m,\] then for every \(k\in \Z\), \[x'=x_0+(b-cm)k,\quad y'=y_0-(b-c)k,\quad z'=z_0+amk,\quad w'=w_0-ak\] is also a solution.
\begin{proof}
The proof is straightforward.
\end{proof}
\end{lem}
\begin{lem}
If \((x_0,y_0, z_0, w_0)\in \Z^4\) is a solution of \[a\sqrt m(x+y\sqrt m)+(b+c\sqrt m)(z+w\sqrt m)=A+B\sqrt m,\] then for every \(l\in \Z\), \[x'=x_0+cml,\quad y'=y_0+bl,\quad z'=z_0-aml,\quad w'=w_0\] is also a solution.
\begin{proof}
The proof is straightforward.
\end{proof}
\end{lem}
\begin{cor}
There is a unique solution \((\bar x,\bar y,\bar z, \bar w)\) of \[a\sqrt m(x+y\sqrt m)+(b+c\sqrt m)(z+w\sqrt m)=A+B\sqrt m\] with \(0\leq \bar z<am\), \(0\leq\bar w<a\).
\begin{proof}
Let  \((x_0,y_0, z_0, w_0)\) be a solution. Let \(k=\lfloor \frac {w_0}{a}\rfloor\). Then \(k\in \Z\), and  \[x'=x_0+(b-cm)k,\quad y'=y_0-(b-c)k,\quad z'=z_0+amk,\quad w'=w_0-ak\] is a solution. By the choice of \(k\), \(0\leq w'<a\) and by Lemma 9, such \(w'\) is unique. \\
Now, let \(l=\lfloor \frac {z'}{am}\rfloor\). Then \(l\in \Z\), and \[\bar x=x'+cml,\quad \bar y=y'+bl,\quad \bar z=z'-aml,\quad \bar w=w'\] is also a solution, and by the choice of \(l\), \(0\leq \bar z<am\). By Lemma 9, such \(\bar z\) is unique. So \((\bar x,\bar y,\bar z, \bar w)\) is a solution with \(0\leq \bar z<am\), \(0\leq\bar w<a\). This is a unique such solution, because if \(\bar z\) and \(\bar w\) are fixed, so are \(\bar x\) and \(\bar y\).
\end{proof}
\end{cor}
\begin{lem} %lemma 12
\(a\sqrt m(x+y\sqrt m)+(b+c\sqrt m)(z+w\sqrt m)=A+B\sqrt m\) has a solution in \(\NM\) if and only if the unique solution \((\bar x,\bar y,\bar z, \bar w)\) with \(0\leq \bar z<am\), \(0\leq\bar w<a\) satisfies \(0\leq \bar x\), \(0\leq \bar y\).
\begin{proof}
Suppose \(\bar x, \bar y\geq 0\). Then \(\bar x+\bar y\sqrt m\), \(\bar z+\bar w\sqrt m\) is a solution in \(\NM\). Now, suppose either \(\bar x<0\) or \(\bar y<0\).  Let \((x', y', z', w')\) be another solution. If \(z'<\bar z\) then \(z'<0\), and if \(w'<\bar w\) then \(w'<0\), by Lemma 9. If \(z'\geq \bar z\) and \(w'\geq \bar w\), then \[x'=\frac{1}{a}(B-cz'-bw')\leq \frac{1}{a}(B-c\bar z-b\bar w)=\bar x,\quad \textrm{and}\]
\[y'=\frac{1}{am}(A-bz'-cmw')\leq \frac{1}{am}(A-b\bar z-cm\bar w)=\bar y\]
so either \(x'<0\) or \(y'<0\).
\end{proof}
\end{lem}
\begin{theorem}
If \(A\geq abm+acm-am-cm-b+1\) and \(B\geq acm+ab-a-b-c+1\), then \(a\sqrt m(x+y\sqrt m)+(b+c\sqrt m)(z+w\sqrt m)=A+B\sqrt m\) has a solution in \(\NM\).
\begin{proof}
Suppose that \(A\geq abm+acm-am-cm-b+1\) and \(B\geq acm+ab-a-b-c+1\). Consider the equation \(a\sqrt m(x+y\sqrt m)+(b+c\sqrt m)(z+w\sqrt m)=A+B\sqrt m\). There is a unique solution \((\bar x,\bar y,\bar z, \bar w)\) with \(0\leq \bar z<am\), \(0\leq\bar w<a\); \((\bar x,\bar y,\bar z, \bar w)\) satisfies 
\[am\bar y+b\bar z+cm\bar w=A,\]
\[a\bar x+c\bar z+b\bar w=B.\] 
Since \(0\leq \bar z\leq am-1\), \(\bar w\leq a-1\), we have 
\[a\bar x=B-c\bar z-b\bar w\geq acm+ab-a-b-c+1-c(am-1)-b(a-1)=1-a\] so
\[\bar x\geq -1+\frac{1}{a}>-1 \iff \bar x\geq 0.\] Likewise, 
\[am\bar y=A-b\bar z-cm\bar w\geq abm+acm-am-cm-b+1-b(am-1)-cm(a-1)\] so
\[\bar y\geq -1+\frac{1}{am}>-1 \iff \bar y\geq 0.\] Therefore \(\bar x+\bar y\sqrt m\in \NM\), \(\bar z+\bar w\sqrt m\in \NM\) is a solution of \(a\sqrt m(x+y\sqrt m)+(b+c\sqrt m)(z+w\sqrt m)=A+B\sqrt m\) in \(\NM\).
\end{proof}
\end{theorem}
\begin{theorem}
\(Frob(a\sqrt m,b+c\sqrt m)=(a\sqrt m-1)(b+c\sqrt m-1)(1+\sqrt m)+\NM\).
\begin{proof}
Note that \[(a\sqrt m-1)(b+c\sqrt m-1)(\sqrt m+1)\]
\[=(abm+acm-am-cm-b+1)+(acm+ab-a-b-c+1)\sqrt m.\]
In view of Theorem 3, to show that \[Frob(a\sqrt m,b+c\sqrt m)=(a\sqrt m-1)(b+c\sqrt m-1)(\sqrt m+1)+\NM\] it suffices to show that there is no \(\beta\in \Z\) such that \[(abm+acm-am-cm-b)+\beta\sqrt m\in Frob(a\sqrt m,b+c\sqrt m)\] and that there is no \(\alpha\in \Z\) such that \[\alpha+(acm+ab-a-b-c)\sqrt m\in Frob(a\sqrt m,b+c\sqrt m).\]
Let \[A=abm+acm-am-cm-b.\] 
Take arbitrary \(k\in \N\), and let \[B=acm+ab-b-c+ak.\] Then \((k,-1,am-1,a-1)\) is the unique solution \((\bar x, \bar y, \bar z, \bar w)\) of \[a\sqrt m(x+y\sqrt m)+(b+c\sqrt m)(z+w\sqrt m)=A+B\sqrt m\] with \(0\leq \bar z<am\), \(0\leq\bar w<a\). By Lemma 12, the equation does not have a solution in \(\NM\). Since the choice of \(k\) was arbitrary, \(B\) can be arbitrarily large. Therefore there is no \(\beta\in \Z\) such that \[(abm+acm-am-cm-b)+\beta\sqrt m\in Frob(a\sqrt m,b+c\sqrt m).\]
Now, let \[B=acm+ab-a-b-c.\]
Take arbitrary \(l\in \N\), and let \[A=abm+acm-cm-b+aml.\]
Then \((-1,l,am-1,a-1)\) is the unique solution \((\bar x, \bar y, \bar z, \bar w)\) of \[a\sqrt m(x+y\sqrt m)+(b+c\sqrt m)(z+w\sqrt m)=A+B\sqrt m\] with \(0\leq \bar z<am\), \(0\leq\bar w<a\). By Lemma 12, the equation does not have a solution in \(\NM\). Since the choice of \(l\) was arbitrary, \(A\) can be arbitrarily large. Therefore there is no \(\alpha\in \Z\) such that \[\alpha+(acm+ab-a-b-c)\sqrt m\in Frob(a\sqrt m,b+c\sqrt m).\]
\end{proof}
\end{theorem}
This completes the solution of the Frobenius problem of the first kind in \(\ZM\) when \(n=2\). We now have the following corollary:
\begin{cor} If \(\alpha_1,\alpha_2\in\NM\setminus\{0\}\), \(\alpha_1,\alpha_2\) span \(1\), and either \(\alpha_i\) has either rational or irrational part \(0\), then \[Frob(\alpha_1,\alpha_2)=(\alpha_1-1)(\alpha_2-1)(1+\sqrt m)+\NM.\]
\end{cor}
\section*{Acknowledgments}
\ \ \ \ Thanks are due to Kari Vaughn, who worked on this problem before the author did and achieved results by brute force that corroborate the general results here.

\end{document}